\pdfoutput=1
\documentclass{amsart}
\usepackage{a4wide}
\usepackage{amsrefs}
\usepackage{amscd}
\usepackage{amsmath}
\usepackage{amssymb}
\usepackage{amsthm}
\usepackage{graphicx}
\usepackage{hyperref}

\theoremstyle{plain}
\newtheorem{theorem}{Theorem}

\newtheorem{fact}{Fact}

\theoremstyle{definition}
\newtheorem{definition}{Definition}

\newtheorem{notation}{Notation}

\theoremstyle{remark}
\newtheorem{remark}{Remark}
\newcommand{\Dom}{\operatorname{Dom}}

\newcommand{\rot}{\operatorname{rot}}
\newcommand{\ind}{\operatorname{ind}}

\begin{document}
\title[Curvature and quantized Arnold strangeness]
{Curvature and quantized Arnold strangeness}
\author{Noboru Ito}
\address{National Institute of Technology, Ibaraki College, 312-8508, Japan}
\email{nito@gm.ibaraki-ct.ac.jp}
\date{November 11, 2022}
\maketitle

\begin{abstract}
By integrating  curvatures multiplied non-trivial densities,  we introduce an integral expression of the \emph{Arnold strangeness} that is a celebrated plane curve invariant.  The key is a partition function by Shumakovitch to reformulate Arnold strangeness.  Our integrating curvatures suggests a quantized Arnold strangeness which Taylor expansion includes the rotation number and the original Arnold strangeness, and also higher terms are invariants of Tabachnikov.      
It is an analogue of the quantization by Viro for Arnold $J^-$ and by  Lanzat-Polyak for $J^+$.    
\end{abstract}
\section{Introduction}\label{intro}
Let $C$ be an oriented generic  immersion $S^1 \to \mathbb{R}^2$, i.e. immersions with a finite set of transversal double points as the only singularities.   
Let $\kappa(t)$ the curvature of $C(t)$.   The following is well-known. 
\begin{fact}[Hopf's Umlaufsatz]\label{fact:hopf}
Let $t \in S^1$ of a plane curve $C : S^1 \to  \mathbb{R}^2$.     
Then the rotation number of $C$ is expressed by 
\[  
\frac{1}{2 \pi} \int_{S^1} \kappa(t) dt .
\]
\end{fact}
\noindent Let $p \in \mathbb{R}^2 \setminus C$ and $\ind_C (p)$ the number of turns made by the vector pointing from $p$ to $C(t)$ ($t \in S^1$) if we follow $C$ along the orientation of $C$, i.e. 
\begin{equation}\label{MappingDegree}
\ind_C (p) := \deg \varphi_p~{\text{where}}~\varphi_p : t \mapsto  \frac{C(t) - p}{|C(t) - p|}.     
\end{equation}
We extend $\ind_C (p)$ to that of  $p \in C$ \cite{Shumakovitch1995, LanzatPolyak2013}.  
If $p$ is a regular point on an edge, $\ind_C (p)$ is the average of values on two regions adjacent to $p$.  
If $p$ is a double point, $\ind_C (p)$ is the average of values on the four regions adjacent to $p$.  This number $\ind_C$ is called the \emph{index} as Viro \cite{Viro1996} does.  
 
The extension of $\ind_C$ to multi-component curves is straightforward.    Therefore, by using the function $\ind_C$, 
Viro \cite{Viro1996} reformulated Arnold $J^-$ \cite{Arnold1994, Arnold1994book} as follows.  Let $\widetilde{C}$ be disjoint circles given by smoothing the double points along the orientation of $C$; let $\Sigma_2$ be the set of regions of $\mathbb{R}^2 \setminus \widetilde{C}$; and let $f_{\sigma}$ be the restriction of $\ind_{\widetilde{C}} |_\sigma$ to $\sigma \in \Sigma_2$.  Clearly, $f_{\sigma}$ is a constant map    
and $\ind_{\widetilde{C}} = \sum_{\sigma} f_{\sigma}$.  Let $\chi(\sigma)$ be the Euler characteristic of $\sigma$ and let $\int_{\mathbb{R}^2 \setminus \widetilde{C}} \ind_{\widetilde{C}} d \chi$ $=$ $\sum_{\sigma} f_{\sigma} (\sigma) \chi(\sigma)$.  
\begin{fact}[{Viro (1996) \cite{Viro1996}}]
\begin{align*} 
{\text{The rotation number of}}~C  &= \int_{\mathbb{R}^2 \setminus \widetilde{C}} \ind_{\widetilde{C}} (x) d \chi (x), \qquad
J^- (C) = 1- \int_{\mathbb{R}^2 \setminus \widetilde{C}} (\ind_{\widetilde{C}} (x))^2 d \chi (x).
\end{align*} 
\end{fact}
\noindent Viro extends the above  formulation to higher terms as  
$
\int_{\mathbb{R}^2 \setminus \widetilde{C}} (\ind_{\widetilde{C}} (x))^r d \chi (x)
$, 
and define their quantized polynomial $P_C (q)$ by  
\[
P_C (q) = \sum_{r=0}^{\infty} \frac{h^r}{r!} \int_{\mathbb{R}^2 \setminus \widetilde{C}} (\ind_{\widetilde{C}} (x))^r d \chi (x)  = \int_{\mathbb{R}^2 \setminus \widetilde{C}} q^{\ind_{\widetilde{C}} (x)} d \chi (x).  
\]
After the work \cite{Viro1996} (1996) was known, Lanzat-Polyak \cite{LanzatPolyak2013} introduced a quantized polynomial $I_q (C)$.  
\begin{fact}[Lanzat-Polyak (2013)  \cite{LanzatPolyak2013}]
For each double point of $C$ $=$ $C(t_1)$ $=$ $C(t_2)$,  letting $\theta_d \in (0, \pi)$ be the non-oriented angle between two tangent vectors $C'(t_1)$ and $-C'(t_2)$, define $I_q (C)$ by 
\begin{equation*}\label{eq:IQ}
I_q (C) = \frac{1}{2 \pi} \left( \int_{\mathbb{S}^1} \kappa (t) \cdot q^{\ind_C} (C(t)) dt  - 
\sum_{d} \theta_d \cdot q^{\ind_{C} (d)} (q^{\frac{1}{2}} - q^{-\frac{1}{2}}) \right).
\end{equation*}
Then its Taylor expansion at $q=1$ satisfies that the first term is the rotation number of $C$ and the second term $I'_1$ satisfies 
$I'_1(C) = \frac{1}{2}(1-J^+(C))$.  
\end{fact}
Since the Arnold basic invariants consist of $J^{+}$, $J^{-}$, and $St$ \cite{Arnold1994book}, it is natural to request an  integral expression by the curvature and its quantization also for $St$.  We discuss this issue here.  

Let $\alpha(t)$ ($t \in C$) and $\alpha_{\pm} (d)$ ($d$ : a double point) be sums of signs defined as in Section~\ref{secPrelim}.       
In the rest of this paper, we suppose that every plane curve $C$ has the base point lies on an exterior edge, which fixes the long curve on $\mathbb{R}^2 \cup \{ \infty \}$ putting the base point on $\{ \infty \}$.   
Let $\rot(C)$ be the rotation number of the corresponding long curve.  
\begin{theorem}\label{mainResult}
Define $St_q (C)$ by     
\[
 \frac{1}{2 \pi (q^{\frac{1}{2}}+q^{-\frac{1}{2}})} \left( \int_{\mathbb{S}^1} \kappa (t)   \alpha(t) q^{ \ind_C (t)} dt  +  
\sum_{d} (\pi - \theta_d) (\alpha_+(d) q^{ \ind_{C} (d) + \frac{1}{2} } - \alpha_-(d) q^{  \ind_{C} (d) - \frac{1}{2}}) 
 \right),   
\]
which becomes a Laurent  polynomial in $\mathbb{Z}[q, q^{-1}]$.  Substituting $q=e^h$ $=$ $\sum_{r=0}^{\infty} \frac{h^r}{r!}$, the $r$th coefficient is an invariant $\frac{St^r(C)}{r!}$ of Tabachnikov.  Taylor expansion of $St_q (C)$ at $q=1$ satisfies that the first term $St_1$ $=$ $\rot(C) $ and the second term $St'_1(C)$ $=$ $St(C)$; $St(C)$ is expressed by       
\begin{equation*}\label{firstSt}
 \frac{1}{4 \pi} \left( \int_{\mathbb{S}^1} \kappa (t)  \alpha(t) \ind_C (t) dt  +  
\sum_{d} (\pi - \theta_d) (\alpha_+(d) \left(\ind_{C} (d) + \frac{1}{2}\right) - \alpha_-(d) \left(\ind_{C} (d) - \frac{1}{2} \right) 
 \right).  
 \end{equation*}
Further,  
by Table~\ref{BehavR}, the comparison is given for three quantizations $P_C (q)$, $I_q (C)$, and $St_q (C)$ corresponding to $J^-$, $J^+$, and $St$, respectively.  
\begin{table}[h!]
\caption{Differences of $P_C (q), I_q (C), St_q (C)$ by positive modifications (Figure~\ref{modification}) 
}\label{BehavR}
\begin{tabular}{|c|c|c|c|}\hline
 & $P_C (q)$ & $I_q (C)$ & $St_q (C)$ \\ \hline
{\rm{direct self-tangency modification}} & $0$ & $-q^{\ind} (q^{\frac{1}{2}} - q^{-\frac{1}{2}})$ & $0$ \\ \hline
{\rm{opposite self-tangency modification}} & $q^{\ind-1} (q-1)^2$ & $0$ & $0$ \\ \hline
{\rm{weak triple-point modification}}  & $0$ & $\frac{1}{2} q^{\ind+\frac{1}{2}} (q^{\frac{1}{2}} - q^{-\frac{1}{2}})^2$ & $q^{\ind}(q-1)$ \\ \hline 
{\rm{strong triple-point modification}} 
& $- q^{\ind-1} (q-1)^3$ & $\frac{1}{2}q^{\ind+\frac{1}{2}} (q^{\frac{1}{2}} - q^{-\frac{1}{2}})^2$ & $q^{\ind}(q-1)$ \\ \hline
\end{tabular}
\end{table}
\end{theorem}
\begin{figure}[htbp] 
\includegraphics[width=12cm]{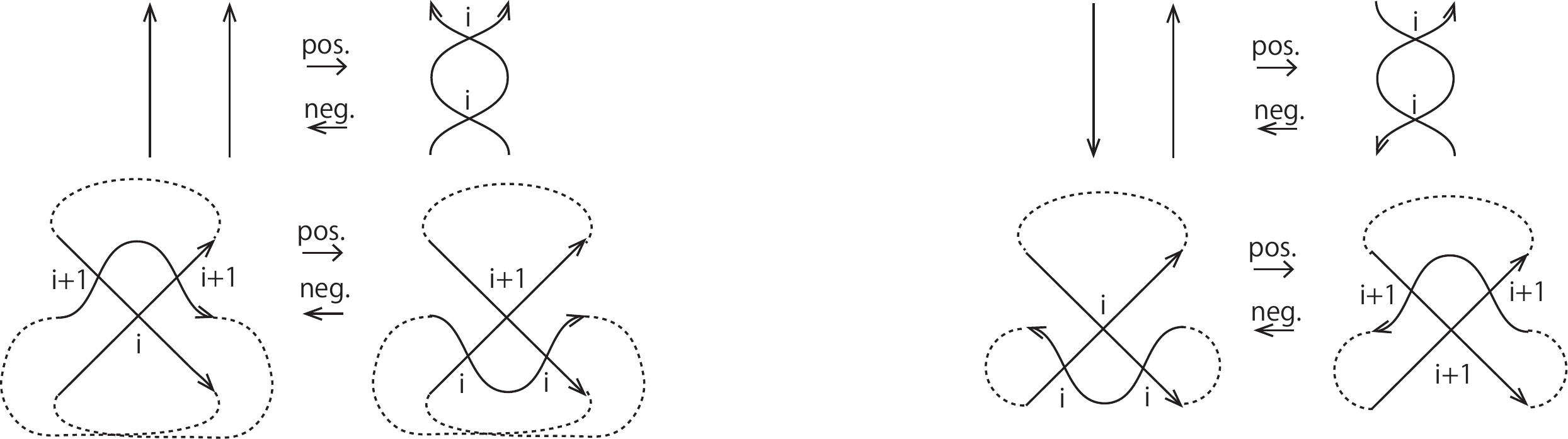} 
   \caption{Modifications with positive/negative directions.  
 Direct/opposite self-tangency modification (upper left/right) and weak/strong triple-point modification (lower left/right) where the case with reversed orientations of the curve is omitted.  Dotted curves indicate connections of curves.}
   \label{modification}
\end{figure}
The locally constant function $\alpha$ of curves is originated in a  \emph{weight} $w$ of Shumakovitch \cite{Shumakovitch1995}, which is related to gleams of Turaev \cite{Turaev1994}.  Using the common weight $w$, 
higher Arnold strangeness $St^r$ is given by Tabachnikov \cite{Tabachnikov1996} for long curves.  Tabachnikov invariant $St^r$ is extended to closed curves by Arakawa-Ozawa \cite {ArakawaOzawa1999}, independently, to closed curves and fronts by Shumakovitch \cite{Shumakovitch1996}.   
\section{A sum $\alpha$ of signs}\label{secPrelim}
\begin{figure}[htbp]
   \includegraphics[width=14cm]{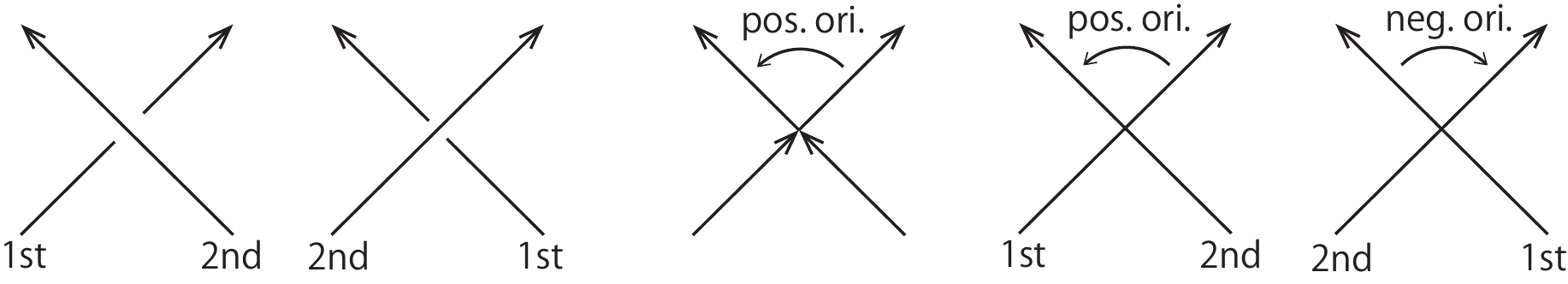} 
\caption{Weights of double points.  Crossings' signs correspond to 
weights of double points (left); the two edges which go into a double point $d$ (center);   
 the double point is negative (positive,~resp.) if the two (the first  and second) edges which go into the double point give a positive (negative,~resp.) orientation (right).    
}
   \label{fig:weight}
\end{figure}

\begin{figure}[htbp]
   \includegraphics[width=14.5cm]{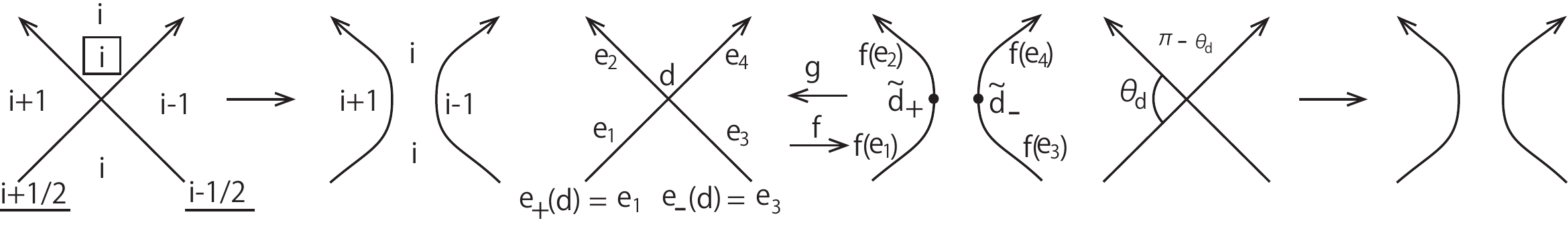}
\caption{Indices of points.  The index $i$ of a double point $d$ is enclosed by the square; indices $i-\frac{1}{2}$, $i+\frac{1}{2}$ are of edges; indices $i-1, i, i+1$ are of regions.   Smoothing a double point and its angle $\theta_d$.}
   \label{fig:angle}
\end{figure}

\begin{definition}[The weight $w$ of a double point $d$ (Shumakovitch  \cite{Shumakovitch1995})]\label{def:weight}
For a plane curve $C$, we choose a base point and we fix an ascending diagram $D_C$ from the base point.   Let $w(d)$ be a local writhe for an ascending diagram.  
\end{definition}
\begin{remark}
Alternatively, the weight is defined by tangent vectors of a  double point.  See Figure~\ref{fig:weight}.  
\end{remark}
\begin{notation}[$\widetilde{d}_+$ and $\widetilde{d}_-$]\label{eLeR}
For a double point $d$, let $e_+ (d)$ and $e_- (d)$ be the two edges
as in Figure~\ref{fig:angle} (center).   If we smooth the double points of $C$, there exists a canonical map $f$  sending each edge $e$ of $C$ to a $1$-simplex of $\widetilde{C}_n$ of an $n$ as in Figure~\ref{fig:angle} (center).   
Then let $\widetilde{d}_+$ ($\widetilde{d}_-$,~resp.) be the vertex where the edge $f(e_+(d))$ ($f(e_-(d))$,~resp.) goes into it.   There is also the canonical map $g$ sending two verticies $\widetilde{d}_+$ and $\widetilde{d}_-$ to the double point $d$.   The symbol $\widetilde{d}_{\epsilon}$ ($\epsilon=+, -$) is often denoted by $\widetilde{d}$ simply.  
\end{notation}
\begin{definition}[$\alpha(\widetilde{C}_n)$; $\alpha(p)$ (a point $p$ on a curve $C(\mathbb{S}^1)$); $\alpha_\pm (d)$ (a double point $d$ on the curve)]\label{def:alpha}
If we smooth the double points of $C$, we have multi-component disjoint simple closed curves $\widetilde{C}$ $=$ $\sqcup_n \widetilde{C}_n$ on the plane.   
Then an integer $\alpha(\widetilde{C}_n)$ is defined by 
\[
\alpha(\widetilde{C}_n) = \rot(\widetilde{C}_n) \sum_{\widetilde{d} \in \widetilde{C}_n} w (\widetilde{d}), 
\]
where $w(\widetilde{d})$ $:=$ $w(d)$ (Definition~\ref{def:weight}) and where $\widetilde{d}$ is as in Notation~\ref{eLeR}.  
Functions $\alpha(p)$ (for any point $p \in C(\mathbb{S}^1)$) and $\alpha_{\pm} (d)$ (for a double point $d \in C(\mathbb{S}^1$)) are defined as follows.     

\noindent $\bullet$ If $p$ is not a double point of $C(\mathbb{S}^1)$, let $\alpha(p)$ $=$ $\alpha(\widetilde{C}_n)$.  

\noindent $\bullet$ If $p$ is a double point $d$, let $\alpha(p)$ ($= \alpha(d)$) $=$ $\frac{1}{2} (\alpha(\widetilde{C}_i) + \alpha(\widetilde{C}_j))$ where $\widetilde{d}_+ \in \widetilde{C}_i$ and $\widetilde{d}_- \in \widetilde{C}_j$.  

\noindent $\bullet$ Let $\alpha_{+} (d)$ $=$ $\alpha(\widetilde{C}_i)$ and $\alpha_{-} (d)$ $=$ $\alpha(\widetilde{C}_j)$ where $\widetilde{d}_+ \in \widetilde{C}_i$ and $\widetilde{d}_- \in \widetilde{C}_j$.   
\end{definition}
\section{Examples of $\alpha$}
\begin{figure}[htbp] 
   \includegraphics[width=10cm]{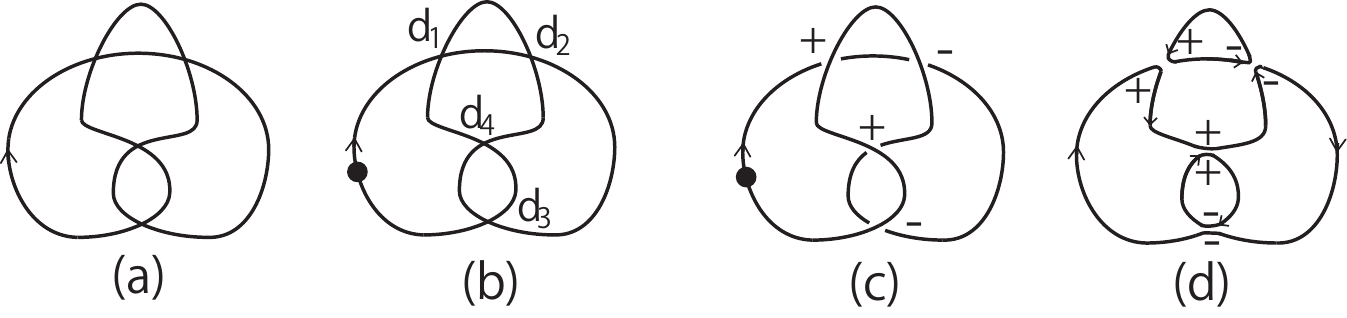} 
   \caption{A computing process to obtain signs, which contributes $ \alpha (\widetilde{C}_n)$ for $\widetilde{C}$ $=$ $\sqcup_n \widetilde{C}_n$.   In this figure, $\alpha (\widetilde{C}_n)$ $=$ $0$ for every $n$.   
   }\label{fig:OneSimplexA}
\end{figure}

\begin{figure}[htbp] 
   \includegraphics[width=10cm]{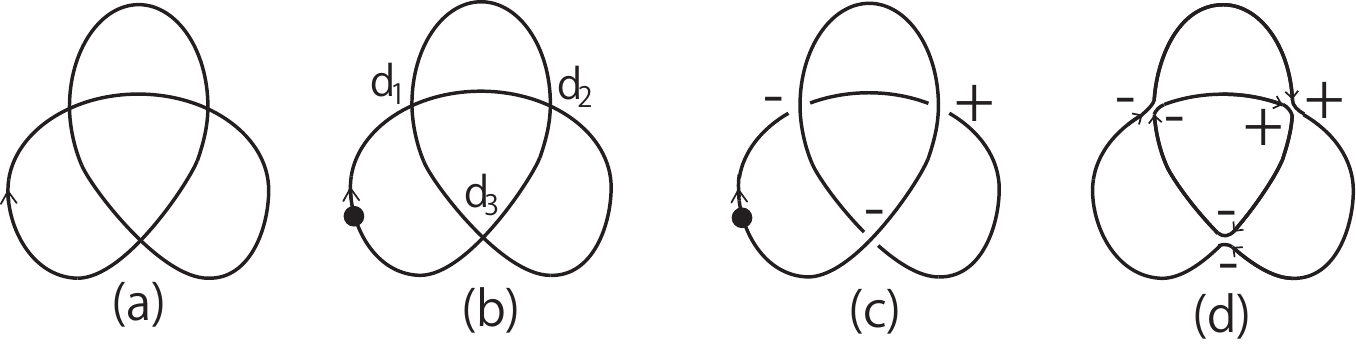} 
   \caption{Another example.   Letting $\widetilde{C}_1$ ($\widetilde{C}_2$,~resp.) be the inner (outer,~resp.) circle, $\rot(\widetilde{C}_n) \alpha(\widetilde{C}_n)$ $=$ $-1+1-1= -1$, i.e. $\alpha(\widetilde{C}_n)$ $=$ $1$, for each $n$.  
}
   \label{fig:OnesimplexTA}
\end{figure}

Figures~\ref{fig:OneSimplexA} and \ref{fig:OnesimplexTA} give examples of computations of $\alpha$.   
\begin{enumerate}
\item[(a)] For a given plane curve $C$ (Figures~\ref{fig:OneSimplexA}, \ref{fig:OnesimplexTA}~(a)), we choose a base point.   
\item[(b)] We assign symbols $d_1, d_2, \dots, d_{2\ell}$ to the double points whose order is derived from a base point (Figures~\ref{fig:OneSimplexA}, \ref{fig:OnesimplexTA}~(b)).  
\item[(c)] We consider the ascending diagram (Figures~\ref{fig:OneSimplexA}), \ref{fig:OnesimplexTA}~(c)) from the base point; it gives the signs of crossings. 
\item[(d)] For each double point, we assign the sign of $d$ to $\widetilde{d}_+$ and $\widetilde{d}_-$.   For $n$, an integer $\alpha(\widetilde{C}_n)$ is given by the sum of signs of $\widetilde{d} \in \widetilde{C}_n$ and $\rot(\widetilde{C}_n)$.  
\end{enumerate}
 
\section{Proof of the formula of the first coefficient in Theorem~\ref{mainResult}}\label{firstProof}
Before we start proving the general case of Theorem~\ref{mainResult}, we prove $St'_1 (C)$ $=$ $St(C)$ since this case is essential.  

Let $\widetilde{C} = \sqcup_n \widetilde{C}_n (\subset \mathbb{R}^2)$, which are disjoint simple closed curves given by smoothing the double points along the orientation of $C$.  Let $U_d$ be a sufficient small neighborhood of a double point $d$ of index $i$ as in Figure~\ref{fig:angle} (left).   
By smoothing $d \in C$, the integral $\int_{S^1} \kappa (t) \alpha(t)  \ind_C (t) dt$ of $\widetilde{C} \cap U_d$ differs from that of $C \cap U_d$ by $\pi - \theta_d$ for the fragment with index $i \pm \frac{1}{2}$ as in Figure~\ref{fig:angle} (right).  Thus this integral part  increases by $(\pi - \theta_d) (\alpha_+(d) ( \ind_{C} (d) + \frac{1}{2} ) - \alpha_-(d) (\ind_{C} (d) - \frac{1}{2}) )$.  
Note that the function $\alpha$ and $\ind_{\widetilde{C}}$ are locally constant; in particular, $ \ind_{\widetilde{C}_n} (t)$ is  constant on $\widetilde{C}_n$; and thus we use the symbol $\ind({\widetilde{C}_n})$ that indicates $\ind_{\widetilde{C}_n} (t)$.   
Let $\Dom(\widetilde{C}_n)$ be the  domain of $\widetilde{C}_n$.  Then

\begin{align*}
&\frac{1}{2 \pi} \left( \int_{\mathbb{S}^1} \kappa (t)  \alpha(t) \ind_C (t) dt  +  
\sum_{d} (\pi - \theta_d) (\alpha_+(d) \left(\ind_{C} (d) + \frac{1}{2}\right) - \alpha_-(d) \left(\ind_{C} (d) - \frac{1}{2} \right) 
 \right) \\
&= \frac{1}{2 \pi} \sum_{n} \int_{\Dom(\widetilde{C}_n)} \kappa (t)  \alpha(t) \ind_{\widetilde{C}} (t) dt \\
&=  \sum_{n} \alpha(\widetilde{C}_n) \ind(\widetilde{C}_n) \frac{1}{2 \pi} \int_{\Dom(\widetilde{C}_n)} \kappa (t) dt \\
\end{align*}
\begin{align*}
&= \sum_{n} \alpha(\widetilde{C}_n) \ind(\widetilde{C}_n) \rot(\widetilde{C}_n) \qquad  (\because {\textrm{Fact}}~\ref{fact:hopf}) \\
&= \sum_{n} \ind(\widetilde{C}_n) \rot^2 (\widetilde{C}_n) \sum_{\widetilde{d} \in \widetilde{C}_n} w(\widetilde{d})   \qquad  (\because \alpha(\widetilde{C}_n) = \rot (\widetilde{C}_n) \sum_{\widetilde{d} \in \widetilde{C}_n} w(\widetilde{d}) \quad ({\rm Definition}~\ref{def:alpha})) \\
&= \sum_{n} \sum_{\widetilde{d} \in \widetilde{C}_n} w(\widetilde{d})  \ind(\widetilde{C}_n) \qquad  (\because \rot^2 (\widetilde{C}_n)=1) \\
&= \sum_{d \in C} w(d)  \left(\ind_C(d)+\frac{1}{2} + \ind_C(d)-\frac{1}{2} \right) \\
&= 2 \sum_{d \in C} w(d) \ind_C(d)  \\
&= 2 St(C).   
\end{align*}
Here we put an explanation for the last two equalities.  The contribution from each double point $d$ into the sum is $w(d) \{ (\ind_C(d)+\frac{1}{2}) + (\ind_C(d)-\frac{1}{2}) \}$ $=$ $2 w(d) \ind_C(d) $, and then we have $2 \sum_{d} w(d) \ind_C(d)$.  The last equality is a known formula by Shumakovitch \cite[Section~1.5]{Shumakovitch1995} because,  throughout of this paper, we suppose that the base point lies on an exterior edge.   
$\hfill \Box$

\section{Proof of the general formula of Theorem~\ref{mainResult}}
Recall Tabachnikov function \cite{Tabachnikov1996}.   
Let $C$ be a generic long curve and $d$ a double point.  For any $k \ge 0$, define $St^k (C)$ by the formula: 
\begin{equation}\label{FTabachnikov}
St^k (C) = \sum_{d} w(d) (\ind_{C} (d))^k , 
\end{equation}
where $w(d)$ is 
as in Definition~\ref{def:weight}.  

Note that every generic long curve  is regarded as a plane curve $C$ with a base point that lies on an exterior edge of $C$ and thus $C$ of $St^k (C)$ is also.  
Using the same notations of Section~\ref{firstProof}, we will show the general formula.  $(q^{\frac{1}{2}}+q^{-\frac{1}{2}}) St_q(C)$ is as follows.  
\begin{align}
\nonumber
&\frac{1}{2 \pi} \left( \int_{\mathbb{S}^1} \kappa (t)  \alpha(t) q^{\ind_C (t)} dt  +  
\sum_{d} (\pi - \theta_d) \left(\alpha_+(d) q^{\ind_{C} (d) + \frac{1}{2}} - \alpha_-(d) q^{\ind_{C} (d) - \frac{1}{2}} \right)
 \right) \\ \nonumber
=& \frac{1}{2 \pi} \sum_{n} \int_{\Dom(\widetilde{C}_n)} \kappa (t)  \alpha(t) q^{\ind_{\widetilde{C}} (t)} dt \\ \nonumber
=& \sum_{n} \alpha(\widetilde{C}_n) q^{\ind(\widetilde{C}_n)} \frac{1}{2 \pi} \int_{\Dom(\widetilde{C}_n)} \kappa (t) dt \\ \nonumber
=& \sum_{n} \alpha(\widetilde{C}_n) q^{\ind(\widetilde{C}_n)} \rot(\widetilde{C}_n) \qquad  (\because {\textrm{Fact}}~\ref{fact:hopf}) \\ \nonumber
=& \sum_{n} q^{\ind(\widetilde{C}_n)} \rot^2 (\widetilde{C}_n) \sum_{\widetilde{d} \in \widetilde{C}_n} w(\widetilde{d})   \qquad  (\because \alpha(\widetilde{C}_n) = \rot (\widetilde{C}_n) \sum_{\widetilde{d} \in \widetilde{C}_n} w(\widetilde{d}) \quad ({\rm Definition}~\ref{def:alpha})) \\ \nonumber
=& \sum_{n} \sum_{\widetilde{d} \in \widetilde{C}_n} w(\widetilde{d})  q^{\ind(\widetilde{C}_n)} \qquad  (\because \rot^2 (\widetilde{C}_n)=1) \\ \label{higherQ}
=&  \sum_{d \in C} w(d) q^{\ind_C(d)} \left(q^{\frac{1}{2}} + q^{-\frac{1}{2}} \right).  
\end{align}
Then  
\[St_q (C) \stackrel{(\ref{higherQ})}{=} \sum_{d \in C} w(d) q^{\ind_C(d)} \stackrel{q=e^h}{=} \sum_{r=0}^{\infty} \frac{h^r}{r!} \sum_{d \in C} w(d) \ind_C^r (d) \stackrel{(\ref{FTabachnikov})}{=} \sum_{r=0}^{\infty} \frac{h^r}{r!} St^r (C). \]  
Further,  by (\ref{higherQ}), substituting $q=1$, $St_1(C)$ $=$ $\sum_{d \in C} w(d)$, which is is the rotation number $\rot(C)$ of the long curve corresponding to $C$.   It is also elementary to see that \[St'_1 (C)= \frac{d(St_q (C))}{dq}|_{q=1} = \sum_{d \in C} w(d) \ind(d)
= St(C),  
\]    
where the last equality is given by  Shumakovitch \cite[Section~1.5]{Shumakovitch1995}.  

Finally, we show Table~\ref{BehavR}.  Since the differences of $P_C (q)$ and $I_q (C)$ are proved in \cite{Viro1996} and \cite{LanzatPolyak2013} respectively, we will check those of $St_q (C)$ for each modification.  Recall that $St_q (C)=\sum_{d \in C} w(d) q^{\ind_C(d)}$ (\ref{higherQ}).  Using Figure~\ref{modification}, direct computation implies the list:    
\begin{itemize}
\item(direct/opposite self-tangency modification)  Let $d, d'$ be two increased double points.  Clearly, $w(d)+w(d')=0$, which implies the difference is $0$.  
\item(weak/strong triple-point modification)  Seeing the weights of three double points in $3 \times 2$ cases corresponding to the possibilities of positions of the base point and curve orientations, it is elementary to check that the difference is  $q^{i+1}-q^{i}$ ($i$ : a value of $\ind_C$ of a double point) in each case.  
\end{itemize}

$\hfill \Box$

\begin{remark}
In \cite{LanzatPolyak2013}, they  compute the case with $\ind=i$ corresponding to Table~\ref{BehavR}.  
\end{remark}
\section*{Acknowledgements}
The author would like to thank Professor~Tomonori Fukunaga for informing the author about \cite{LanzatPolyak2013}.   The author also thank Professors Mitsuhiro Imada and Yumiko Takei for fruitful discussions.  
The work was partially supported by JSPS KAKENHI Grant Numbers JP20K03604 and JP22K03603.  

\bibliographystyle{plain}
\bibliography{RefA}
\end{document}